\documentstyle[12pt]{article}
\newcommand{\nc}{\newcommand}
\nc{\al}{\alpha}
\nc{\g}{\gamma}
\nc{\G}{\Gamma}
\nc{\D}{\Delta}
\nc{\la}{\lambda}
\nc{\La}{\Lambda}
\nc{\var}{\varphi}
\nc{\kvt}{\sqrt{t}}
\nc{\hn}{h^\vee}
\nc{\kn}{k^\vee}
\nc{\pa}{\partial}
\nc{\nn}{\nonumber \\ }
\nc{\hf}{\frac{1}{2}}         
\nc{\fabc}{{f_{ab}}^c}
\nc{\binomial}[2]{\left (\begin{array}{c} {#1}\\ {#2} \end{array}
\right )}
\nc{\ben}{\begin{equation}}
\nc{\een}{\end{equation}}
\nc{\bea}{\begin{eqnarray}}
\nc{\eea}{\end{eqnarray}}
\nc{\bra}[1]{\langle {#1}|}
\nc{\ket}[1]{|{#1}\rangle}
\newcommand{\Z}{\mbox{$Z\hspace{-2mm}Z$}}
\nc{\C}{\mbox{\hspace{1.24mm}\rule{0.2mm}{2.5mm}\hspace{-2.7mm} C}}
\nc{\Nat}{\mbox{\hspace{.04mm}\rule{0.2mm}{2.8mm}\hspace{-1.5mm} N}}
\newcommand{\R}{\mbox{\hspace{.04mm}\rule{0.2mm}{2.8mm}\hspace{-1.5mm} R}}

\nc{\spa}{\hspace{1 cm},\hspace{1 cm}}
\nc{\vs}{\vspace}
\nc{\NP}[1]{Nucl.\ Phys.\ {\bf #1}}
\nc{\PL}[1]{Phys.\ Lett.\ {\bf #1}}
\nc{\CMP}[1]{Commun.\ Math.\ Phys.\ {\bf #1}}
\nc{\PR}[1]{Phys.\ Rev.\ {\bf #1}}
\nc{\PRL}[1]{Phys.\ Rev.\ Lett.\ {\bf #1}}
\nc{\PTP}[1]{Prog.\ Theor.\ Phys.\ {\bf #1}}
\nc{\PTPS}[1]{Prog.\ Theor.\ Phys.\ Suppl.\ {\bf #1}}
\nc{\MPL}[1]{Mod.\ Phys.\ Lett.\ {\bf #1}}
\nc{\IJMP}[1]{Int.\ Jour.\ Mod.\ Phys.\ {\bf #1}}
\nc{\IM}[1]{Invent.\ Math.\ {\bf #1}}
\nc{\SJNP}[1]{Sov. J. Nucl. Phys.\ {\bf #1}}

\begin{document}

\topmargin -5mm
\oddsidemargin 5mm

\begin{titlepage}
\setcounter{page}{0}
\begin{flushright}
NBI-HE-97-31\\
July 1997
\end{flushright}

\vs{8mm}
\begin{center}
{\Large Explicit Decompositions of Weyl Reflections}\\[.3cm]
{\Large in Affine Lie Algebras}

\vs{8mm}
{\large J{\o}rgen Rasmussen}\footnote{e-mail address: 
jrasmussen@nbi.dk}\\[.2cm]
{\em The Niels Bohr Institute, Blegdamsvej 17, DK-2100 Copenhagen \O,
Denmark}\\[.5cm]

\end{center}

\vs{8mm}
\centerline{{\bf{Abstract}}}
\noindent In this paper explicit decompositions are provided of 
the Weyl reflections in affine Lie algebras,
in terms of fundamental Weyl reflections.\\[.4cm]

\end{titlepage}
\newpage
\renewcommand{\thefootnote}{\arabic{footnote}}
\setcounter{footnote}{0}

\section{Introduction}

An understanding of the Weyl group of an affine Lie algebra resides
in the basis of affine Lie algebra theory \cite{Kac}. 
Just as in the case of the usual
finite dimensional Lie algebras, the (affine)
Weyl group is fundamental in discussions
on e.g. characters. In applications of the Weyl group it is sometimes of 
importance to be able to decompose the Weyl reflections 
into products of fundamental reflections which are reflections with
respect to simple roots. This is the case when considering singular
vectors along the lines of Malikov, Feigin and Fuks \cite{MFF}.

The main result in this paper is the presentation of explicit 
decompositions into fundamental reflections of all Weyl reflections in 
all affine Lie algebras based on simple finite dimensional Lie algebras
of the types $A_r,\ B_r,\ C_r,\ D_r,\ E_6,\ E_7,\ E_8,\ F_4$ and $G_2$. 
The decompositions we present
rely on a new universal lemma and on well known explicit,
algebra dependent realizations of the root systems in the associated
finite dimensional Lie algebras (see e.g. \cite{Fuc}).
The lemma reduces
the problem to an associated problem of decomposing classical
Weyl reflections. By classical Weyl reflections we mean reflections in a 
standard finite dimensional Lie algebra. Explicit decompositions of
the classical Weyl reflections are also worked out.
To the best of our knowledge, explicit decompositions of affine Weyl
reflections are only known in a few examples, see e.g. \cite{Dob, DS}.

Our motivation for considering explicit decompositions of affine Weyl
reflections is the wish to use them in a study of fusion rules in conformal
field theories based on affine current algebras \cite{Ras3}. In that respect
explicit decompositions are necessary in order to generalize a work
by Awata and Yamada \cite{AY} on fusion rules for affine $SL(2)$ current
algebra. Recently, their approach has been used to determine the fusion rules
for admissible representations \cite{KK,KW} of
affine $OSp(1|2)$ current algebra \cite{ER}. Moreover, one needs differential 
operator realizations of the underlying finite dimensional Lie algebras.
In the work \cite{Ras1,PRY} by Petersen, Yu and the present author, such
realizations have been worked out for all simple finite dimensional Lie
algebras. The decompositions presented in this paper then allow generalizing
the work \cite{AY}. 

The remaining part of the
paper is organized as follows. Section 2 serves to fix notation. Some  
basic Lie algebra and affine Lie algebra properties are reviewed. 
In Section 3 explicit decompositions of affine Weyl reflections are 
worked out first in terms of simple reflections and classical
Weyl reflections of the
associated finite dimensional Lie algebra. The classical Weyl reflections are
then decomposed into fundamental reflections.
Section 4 contains concluding
remarks, while explicit, algebra dependent realizations of root systems in
finite dimensional Lie algebras are reviewed in Appendix A.  

\section{Notation}

\subsection{Lie Algebras}
Let {\bf g} be a simple Lie algebra of rank {\bf g} = $r$.
{\bf h} is a Cartan subalgebra of {\bf g}. The set of (positive) roots
is denoted ($\Delta_+$) $\Delta$, and we write $\al>0$ if $\al\in\Delta_+$.
The simple roots are $\{\al_i\}_{i=1,...,r}$. $\theta$ is the highest root,
while $\al^\vee = 2\al/\al^2$ is the root dual to $\al$. 
Using the triangular decomposition 
\ben
 \mbox{{\bf g}}=\mbox{{\bf g}}_-\oplus\mbox{{\bf h}}\oplus\mbox{{\bf g}}_+
\een
the raising and lowering generators are denoted $e_\al\in$ {\bf g}$_+$ and
$f_\al\in$ {\bf g}$_-$ respectively with $\al\in\Delta_+$, and 
$h_i\in$ {\bf h} are the Cartan generators. 
We let $j_a$ denote an arbitrary Lie algebra generator. 
For simple roots one sometimes writes $e_i=e_{\al_i}, f_i=f_{\al_i}$.
The 3$r$ generators $e_i,h_i,f_i$ are the Chevalley generators.
Their commutator relations are
\bea
\left[h_i,h_j\right]=0&&\left[e_i,f_j\right]=\delta_{ij}h_j\nn
  \left[h_i,e_j\right]=A_{ij}e_j&&\left[h_i,f_j\right]=-A_{ij}f_j
\eea
where $A_{ij}$ is the Cartan matrix. 
In the Cartan-Weyl basis we have
\ben
 [h_i,e_\al]=(\al_i^\vee,\al)e_\al\spa [h_i,f_\al]=
-(\al_i^\vee,\al)f_\al
\een
and
\ben
 \left[e_\al,f_\al\right]=h_\al=G^{ij}(\al_i^\vee,\al^\vee)h_j
\een
where the metric $G_{ij}$ is related to the Cartan matrix as
$A_{ij}=\al_i^\vee\cdot\al_j=(\al_i^\vee,\al_j)=
G_{ij}\al_j^2/2$, while the non-vanishing elements of the
Cartan-Killing form are
\ben
  \kappa_{\al,-\beta}=\kappa(e_\al f_{\beta})=\frac{2}{\al^2}
  \delta_{\al,\beta}\spa \kappa_{ij}=\kappa(h_i h_j)=G_{ij}
\een

The Weyl reflections $\sigma_\al$ acting on $\la\in\mbox{{\bf h}}^*$ 
(or equivalently on $\beta\in\D$) are defined by
\ben
 \sigma_\al(\la)=\la-(\la,\al^\vee)\al
\een
and are generators of the Weyl group. It is easily seen that 
\ben
 \sigma_\al^{-1}=\sigma_\al\spa \sigma_{-\al}=\sigma_\al\spa
  \sigma_{\sigma_\al(\beta)}=\sigma_\al\circ\sigma_\beta\circ\sigma_\al
\label{easy}
\een

\subsection{Affine Lie Algebras}

The loop algebra {\bf g}$_t=$ {\bf g} $\otimes\ \C[t,t^{-1}]$ is generated
by the elements $j_a(n)=j_a\otimes t^n$ with defining commutator 
relations
\ben 
 [j_a(n),j_b(m)]=[j_a,j_b](n+m)
\een
It may be centrally extended by adding to it a central algebra element
$k$ which may be treated as a constant since it commutes with all 
generators. Denoting the generators by capital letters we have the following
commutator relations of the extended algebra
\ben
 [J_a(n),J_b(m)]=\fabc J_c(n+m)+\kappa_{ab}kn\delta_{n+m,0}
\een
By further including a derivation $D$ satisfying
\ben
 [D,J_a(n)]=nJ_a(n)\spa [D,k]=0
\een
one has obtained an affine Lie algebra of level $\kn$ where
\ben
 \kn=\frac{2k}{\theta^2}
\een
Note that whenever we consider affine Lie algebras in connection with the
Virasoro algebra (through the Sugawara construction) we may take $D=-L_0$.
Also note that the generators $\{J_a(n)\}$ are the modes of the currents
$\bar{J}_a(z)$ in an affine current algebra
\bea
 \bar{J}_a(z)\bar{J}_b(w)&=&\frac{\kappa_{ab}k}{(z-w)^2}+\frac{\fabc 
  \bar{J}_c(w)}{z-w}+
 \mbox{regular\ terms}\nn
 \bar{J}_a(z)&=&\sum_{n\in\Z}J_a(n)z^{-n-1}
\eea

The roots $\hat{\al}$ with respect to $(H,k,D)$ are 
\ben
 \hat{\al}=\hat{\al}(n)=(\al,0,n)
\label{realroot}
\een
and 
\ben
 \hat{\al}=n\delta=(0,0,n)\spa n\in\Z\setminus\{0\}
\een
The simple roots are 
\bea
 \hat{\al}_i&=&(\al_i,0,0)\spa 1\leq i\leq r\nn
 \hat{\al}_0&=&(-\theta,0,1)
\eea
and thus, the positive roots are
\ben
 \hat{\al}=(\al,0,n)\spa n>0\ \ \mbox{or}\ \ (n=0,\ \al>0)
\een

Weyl reflections (with respect to the roots (\ref{realroot})) of the weight 
lattice of an affine Lie algebra are defined by
\ben
 \sigma_{\hat{\al}}(\La)=\La-\La\cdot\hat{\al}^\vee\hat{\al}
\een
where the scalar product of 2 elements $\La=(\la,k,m)$ and
$\La'=(\la',k',m')$ is defined by
\ben
 (\La,\La')=\La\cdot\La'=\la\cdot\la'+km'+k'm
\een
Reflections with respect to the simple roots are sometimes written
$\sigma_i$, $i=0,...,r$. It is obvious that $\sigma_{\hat{\al}}^{-1}=\sigma_{
\hat{\al}}$, and it is easily seen that the Weyl reflection with respect
to $\hat{\al}=(\al,0,n)$ takes the weight $\La=(\la,k,m)$ to
\bea
 \sigma_{\hat{\al}}(\La)&=&\left(\sigma_\al(\la_n(\al^\vee)),k,m+\frac{1}{2k}
  (\la^2-\la_n^2(\al^\vee))\right)\nn
 \la_n(\al^\vee)&=&\la+nk\al^\vee
\eea
Inspired by this one may introduce the translation operators
$t_{\al^\vee}$ defined by
\ben
 t_{\al^\vee}(\La)=\left(\la_1(\al^\vee),k,m+\frac{1}{2k}(\la^2-
  \la_1^2(\al^\vee))\right)
\een
By induction it is seen that for all $n\in \Z$
\ben
 t_{\al^\vee}^n(\La)=\left(\la_n(\al^\vee),k,m+\frac{1}{2k}(\la^2-
  \la_n^2(\al^\vee))\right)
\een
The translation operators satisfy
\ben
 t_{\al^\vee+\beta^\vee}=t_{\al^\vee}\circ t_{\beta^\vee}=t_{\beta^\vee}
  \circ t_{\al^\vee}\spa t_{-\al^\vee}
  =t_{\al^\vee}^{-1}
\een
and in particular
\ben
 \sigma_{\hat{\al}}=\sigma_\al\circ t_{\al^\vee}^n
\label{split}
\een
where the definition of the Weyl reflection $\sigma_\al$ is trivially
generalized to the following action on $r+2$ dimensional weights
$\La=(\la,k,m)$:
\ben
 \sigma_\al(\La)=\sigma_{(\al,0,0)}(\La)=(\sigma_\al(\la),k,m)
\een
Furthermore, the translation operator satisfies 
\ben
 t_{\sigma_\al(\beta^\vee)}=\sigma_\al\circ t_{\beta^\vee}\circ\sigma_\al
\een
and from this it follows immediately that
\ben
 t_{\al^\vee}=\sigma_{(-\al,0,1)}\circ\sigma_{(\al,0,0)}
\een
and in particular
\ben
 t_{\theta^\vee}=\sigma_0\circ\sigma_\theta
\een
Thus, in the case of $SL(2)$ where $r=1$ and hence $\theta=\al_1$, a general
(positive) affine Weyl reflection may be decomposed as
\ben
 \sigma_{(\theta,0,n)}=\sigma_1\circ(\sigma_0\circ\sigma_1)^n
\een

\section{Decompositions of Affine Weyl Reflections}

Let us first present the following lemma which is a key ingredient in our
decomposition procedure of affine Weyl reflections\\[.2cm]
{\bf Lemma}
\ben
 t_{\al^\vee}=\sigma_0\circ\sigma_{\theta-\al}\circ\sigma_0
  \circ\sigma_\al\spa \theta^\vee\cdot\al=1
\een
Even though the proof of the lemma is straightforward, it seems to be a new 
result. Note that the condition $\theta^\vee\cdot\al=1$ implies
$\theta-\al>0$ (see e.g. \cite{MP}).

A general affine Weyl reflection $\sigma_{\hat{\al}}$, $\hat{\al}=(\al,0,n)$,
may be written 
\ben
 \sigma_{\hat{\al}}=\sigma_\al\circ t_{\al^\vee}^n=\sigma_\al\circ
  \left(t_{\al_{i_1}^\vee}\circ...\circ t_{\al_{i_N}^\vee}\right)^n
\een
where
\ben
 \al^\vee=\al_{i_1}^\vee+...+\al_{i_N}^\vee
\een
is a unique (integer) expansion of the dual root $\al^\vee$ on dual simple 
roots $\al_i^\vee$. Here we have assumed that $\al>0$. For $\hat{\al}=
(-\al,0,n)$ with $\al>0$, we have
\ben
 \sigma_{\hat{\al}}=\sigma_{-\al}\circ t_{-\al^\vee}^n=\sigma_\al\circ
  \left(t_{\al_{i_1}^\vee}\circ...\circ t_{\al_{i_N}^\vee}\right)^{-n}
\een
This reduces our first problem of decomposing the affine Weyl reflections
in terms of fundamental reflections and classical Weyl reflections,
to one of determining the translation operators
with respect to dual simple roots.

Using the explicit representations of the root systems in Appendix A, it
turns out that in the (generally non-unique) representation 
\ben
 \al_i^\vee=\theta^\vee-\beta_1^{(i)\vee}-...-\beta_{M_i}^{(i)\vee}
\een
of $\al_i^\vee$, one needs at most 2 dual positive roots $\beta_j^{(i)\vee}$,
that is $M=1,2$ ($M=0$ only in the simple case of $A_1\simeq sl(2)$). This
means that
\ben
 t_{\al_i^\vee}=t_{\theta^\vee-\beta_1^{(i)\vee}-\beta_2^{(i)\vee}}=
  t_{\beta_1^{(i)\vee}}^{-1}\circ 
  t_{\beta_2^{(i)\vee}}^{-1}\circ t_{\theta^\vee}
\een
Example by example one may now attempt to use the lemma on $t_{\beta_1^\vee}$
(and for $M=2$ also on $t_{\beta_2^\vee}$). It turns out that in most cases
the lemma applies. However, for $t_{\al_i^\vee}$, $1<i<r$, in $C_r$ (and
for $t_{\al_1^\vee}$ in $C_r$ and for $t_{\al_2^\vee}$ in $G_2$, see below)
it does not apply directly, but one may then repeat the 
procedure for $\beta_j^\vee=\theta^\vee-\g_1^\vee-...-\g_{M'}^\vee$.
Indeed we may write
\bea
 \al_i^\vee&=&\theta^\vee-(\al_{1i}^-)^\vee-\al_{i+1,i+1}^\vee\nn
 t_{\al_i^\vee}&=&t_{(\al_{1i}^-)^\vee}^{-1}\circ t_{\al_{i+1,i+1}^\vee}^{-1}
  \circ t_{\theta^\vee}
\eea
where $\theta^\vee\cdot\al_{1i}^-=1$, but 
$\theta^\vee\cdot\al_{i+1,i+1}\neq1$. In the second step we write
\bea
 \al_{i+1,i+1}^\vee&=&\theta^\vee-(\al_{1,i+1}^-)^\vee\spa\theta^\vee\cdot
  \al_{1,i+1}^-=1\nn
 t_{\al_{i+1,i+1}^\vee}&=&t_{\theta^\vee}\circ t_{(\al_{1,i+1}^-)^\vee}^{-1}
\eea
and in conclusion we have
\bea
 t_{\al_i^\vee}&=&t_{(\al_{1i}^-)^\vee}^{-1}\circ t_{(\al_{1,i+1}^-)^\vee}\nn
 &=&\sigma_{\al_{1i}^-}\circ\sigma_0\circ\sigma_{\al_{1i}^+}
  \circ\sigma_{\al_{1,i+1}^+}\circ\sigma_0\circ\sigma_{\al_{1,i+1}^-}
\eea
In the following Table 1 and Table 2 we have summarized 
our findings.
\newpage
\noindent {\bf Table 1}
\begin{center}
\begin{tabular}{|c||lr|}\hline
 &&\\
 {\bf g}&Fundamental translations  &\\ &&\\ \hline\hline
 &&\\
 $A_r$ &
  $t_{\al_1^\vee}=\sigma_{\al_{2,r+1}}\circ\sigma_0\circ
  \sigma_{\al_1}\circ\sigma_\theta$ &\\
 & $t_{\al_i^\vee}=\sigma_{\al_{1i}}\circ\sigma_0\circ\sigma_{\al_{i,r+1}}
  \circ\sigma_0\circ\sigma_{\al_{i+1,r+1}}\circ\sigma_0\circ
  \sigma_{\al_{1,i+1}}
 \circ\sigma_\theta$  & $1<i<r$\\ 
 & $t_{\al_r^\vee}=\sigma_{\al_{1r}}\circ\sigma_0\circ\sigma_{\al_r}\circ
  \sigma_\theta     $ &\\
 &&\\ \hline
 &&\\
  $B_r$ &$t_{\al_1^\vee}=\sigma_{\al_{22}}\circ\sigma_0\circ\sigma_{\al_{11}}
  \circ\sigma_\theta$  &\\ 
 & $t_{\al_2^\vee}=\sigma_{\al_{13}^+}\circ\sigma_0\circ\sigma_{\al_2}
  \circ\sigma_\theta$  &\\
 & $t_{\al_i^\vee}=\sigma_{\al_{1i}^-}\circ\sigma_0\circ\sigma_{\al_{2i}^+}
  \circ\sigma_0\circ\sigma_{\al_{2,i+1}^+}\circ\sigma_0\circ\sigma_{\al_{
  1,i+1}^-}\circ\sigma_\theta$  &  $2<i<r$\\
 & $t_{\al_r^\vee}=\sigma_{\al_{1r}^-}\circ\sigma_0\circ\sigma_{\al_{2r}^+}
  \circ\sigma_0\circ\sigma_{\al_{2r}^-}\circ
  \sigma_0\circ\sigma_{\al_{1r}^+}\circ\sigma_\theta$  &\\ 
 &&\\ \hline &&\\
  $C_r$ &$t_{\al_1^\vee}=\sigma_0\circ\sigma_{\al_{12}^+}\circ\sigma_0
  \circ\sigma_{\al_1}          $  &\\ 
 & $t_{\al_i^\vee}=\sigma_{\al_{1i}^-}\circ\sigma_0\circ\sigma_{\al_{1i}^+}
  \circ\sigma_{\al_{1,i+1}^+}\circ\sigma_0\circ\sigma_{\al_{1,i+1}^-}$
  &$1<i<r$\\
 & $t_{\al_r^\vee}=\sigma_{\al_{1r}^-}\circ\sigma_0\circ\sigma_{\al_{1r}^+}
  \circ\sigma_\theta$&\\
 &&\\ \hline
 &&\\
  $D_r$ &$t_{\al_1^\vee}=\sigma_{\al_{2r}^-}\circ\sigma_0\circ
  \sigma_{\al_{1r}^+}\circ\sigma_0\circ\sigma_{\al_{2r}^+}\circ\sigma_0\circ
  \sigma_{\al_{1r}^-}\circ\sigma_\theta          $  &\\ 
 & $t_{\al_2^\vee}=\sigma_{\al_{13}^+}\circ\sigma_0\circ\sigma_{\al_2}
  \circ\sigma_\theta$ &\\
 & $t_{\al_i^\vee}=\sigma_{\al_{1i}^-}\circ\sigma_0\circ\sigma_{\al_{2i}^+}
  \circ\sigma_0\circ\sigma_{\al_{2,i+1}^+}\circ\sigma_0\circ\sigma_{
  \al_{1,i+1}^-}\circ\sigma_\theta$ & $2<i<r$\\
 & $t_{\al_r^\vee}=\sigma_{\al_{1,r-1}^-}\circ\sigma_0\circ\sigma_{
  \al_{2,r-1}^+}\circ\sigma_0\circ\sigma_{\al_{2r}^-}\circ\sigma_0
  \circ\sigma_{\al_{1r}^+}\circ\sigma_\theta$ &\\
 &&\\ \hline
\end{tabular}
\end{center}
\newpage
\noindent {\bf Table 2}
\begin{center}
\begin{tabular}{|c||lr|}\hline
 &&\\
 {\bf g}&Fundamental translations  &\\ &&\\ \hline\hline
 &&\\
  $E_6$ &$t_{\al_1^\vee}=\sigma_{\al_{+--++}}\circ\sigma_0\circ
  \sigma_{\al_{23}^+}\circ\sigma_0\circ\sigma_{\al_{13}^+}\circ\sigma_0\circ
  \sigma_{\al_{-+-++}}\circ\sigma_\theta          $  &\\ 
 & $t_{\al_2^\vee}=\sigma_{\al_{++--+}}\circ\sigma_0\circ
  \sigma_{\al_{34}^+}\circ\sigma_0\circ\sigma_{\al_{24}^+}\circ\sigma_0\circ
  \sigma_{\al_{+-+-+}}\circ\sigma_\theta          $  &\\ 
 & $t_{\al_3^\vee}=\sigma_{\al_{+++--}}\circ\sigma_0\circ
  \sigma_{\al_{45}^+}\circ\sigma_0\circ\sigma_{\al_{35}^+}\circ\sigma_0\circ
  \sigma_{\al_{++-+-}}\circ\sigma_\theta          $  &\\ 
 & $t_{\al_4^\vee}=\sigma_{\al_{-+++-}}\circ\sigma_0\circ
  \sigma_{\al_{15}^+}\circ\sigma_0\circ\sigma_{\al_{14}^+}\circ\sigma_0\circ
  \sigma_{\al_{-++-+}}\circ\sigma_\theta          $  &\\ 
 & $t_{\al_5^\vee}=\sigma_{\al_{23}^+}\circ\sigma_0\circ\sigma_{\al_{+--++}}
  \circ\sigma_0\circ\sigma_{\al_{45}^+}\circ\sigma_0\circ
  \sigma_{\al_{+++--}}\circ\sigma_\theta          $  &\\ 
 & $t_{\al_6^\vee}=\sigma_{\al_{--+++}}\circ\sigma_0\circ\sigma_{\al_6}\circ
  \sigma_\theta$ &\\
 &&\\ \hline
 &&\\
 $E_7$ &$t_{\al_1^\vee}=\sigma_{\al_{+-++++}}\circ\sigma_0\circ
  \sigma_{\al_{-+----}}\circ\sigma_0\circ
  \sigma_{\al_{+-----}}\circ\sigma_0\circ
  \sigma_{\al_{-+++++}}\circ\sigma_\theta          $  &\\ 
 &$t_{\al_2^\vee}=\sigma_{\al_{++-+++}}\circ\sigma_0\circ
  \sigma_{\al_{--+---}}\circ\sigma_0\circ
  \sigma_{\al_{-+----}}\circ\sigma_0\circ
  \sigma_{\al_{+-++++}}\circ\sigma_\theta          $  &\\ 
 &$t_{\al_3^\vee}=\sigma_{\al_{+++-++}}\circ\sigma_0\circ
  \sigma_{\al_{---+--}}\circ\sigma_0\circ
  \sigma_{\al_{--+---}}\circ\sigma_0\circ
  \sigma_{\al_{++-+++}}\circ\sigma_\theta          $  &\\ 
 &$t_{\al_4^\vee}=\sigma_{\al_{++++-+}}\circ\sigma_0\circ
  \sigma_{\al_{----+-}}\circ\sigma_0\circ
  \sigma_{\al_{---+--}}\circ\sigma_0\circ
  \sigma_{\al_{+++-++}}\circ\sigma_\theta          $  &\\ 
 &$t_{\al_5^\vee}=\sigma_{\al_{+++++-}}\circ\sigma_0\circ
  \sigma_{\al_{-----+}}\circ\sigma_0\circ
  \sigma_{\al_{----+-}}\circ\sigma_0\circ
  \sigma_{\al_{++++-+}}\circ\sigma_\theta          $  &\\ 
 &$t_{\al_6^\vee}=\sigma_{\al_{-+++++}}\circ\sigma_0\circ
  \sigma_{\al_6}\circ\sigma_\theta          $  &\\ 
 &$t_{\al_7^\vee}=\sigma_{\al_{--+++-}}\circ\sigma_0\circ
  \sigma_{\al_{++---+}}\circ\sigma_0\circ
  \sigma_{\al_{-----+}}\circ\sigma_0\circ
  \sigma_{\al_{+++++-}}\circ\sigma_\theta          $  &\\ 
 &&\\ \hline
 &&\\
  $E_8$ &$t_{\al_i^\vee}=\sigma_{\al_{i7}^+}\circ\sigma_0\circ\sigma_{
  \al_{i8}^-}\circ\sigma_0\circ\sigma_{\al_{i+1,8}^-}\circ\sigma_0\circ
  \sigma_{\al_{i+1,7}^+}\circ\sigma_\theta$  &$1\leq i\leq5$\\ 
  &$t_{\al_6^\vee}=\sigma_{\al_{68}^+}\circ\sigma_0\circ\sigma_{\al_6}
  \circ\sigma_\theta$  &\\
  &$t_{\al_7^\vee}=\sigma_{\al_{--+++++}}\circ\sigma_0\circ\sigma_{
  \al_{++----+}}\circ\sigma_0\circ\sigma_{\al_{27}^+}\circ\sigma_0\circ
  \sigma_{\al_{28}^-}\circ\sigma_\theta$ &\\
  &$t_{\al_8^\vee}=\sigma_{\al_{17}^-}\circ\sigma_0\circ\sigma_{\al_{18}^+}
  \circ\sigma_0\circ\sigma_{\al_{28}^-}\circ\sigma_0\circ\sigma_{\al_{27}^+}
  \circ\sigma_\theta$ &\\  
 &&\\ \hline
 &&\\
  $F_4$ &$t_{\al_1^\vee}=\sigma_{\al_{23}^+}\circ\sigma_0\circ\sigma_{
  \al_{13}^-}\circ\sigma_0\circ\sigma_{\al_{24}^+}\circ\sigma_0
  \circ\sigma_{\al_{14}^-}\circ\sigma_\theta$ &\\
 & $t_{\al_2^\vee}=\sigma_{\al_{13}^+}\circ\sigma_0\circ\sigma_{\al_2}
  \circ\sigma_\theta$  &\\
 & $t_{\al_3^\vee}=\sigma_{\al_{13}^-}\circ\sigma_0\circ\sigma_{\al_{23}^+}
  \circ\sigma_0\circ\sigma_{\al_{24}^+}\circ\sigma_0\circ\sigma_{\al_{14}^-}
  \circ\sigma_\theta$ &\\
 & $t_{\al_4^\vee}=\sigma_{\al_{14}^-}\circ\sigma_0\circ\sigma_{\al_{24}^+}
  \circ\sigma_0\circ\sigma_{\al_{24}^-}\circ\sigma_0\circ\sigma_{\al_{14}^+}
  \circ\sigma_\theta$ &\\
 &&\\ \hline
 &&\\
  $G_2$ &$t_{\al_1^\vee}=\sigma_{\al_2}\circ\sigma_0\circ\sigma_{\al_{+-+}}
  \circ\sigma_0\circ\sigma_{\al_2}\circ\sigma_0\circ\sigma_{\al_{+-+}}
  \circ\sigma_\theta         $  &\\ 
 & $t_{\al_2^\vee}=\sigma_0\circ\sigma_{\al_{+-+}}\circ\sigma_0
  \circ\sigma_{\al_2}$ &\\
 &&\\ \hline
\end{tabular}
\end{center}
\newpage
\noindent
Note that these expressions are neither unique, nor necessarily irreducible. 
First of all because our construction itself is not unique. Let us illustrate 
this by considering some examples. A simple one is provided by 
$t_{\al_1^\vee}$ in $E_7$ where $\al_1^\vee$ may 
be represented in several ways
\bea
 \al_1^\vee&=&\theta^\vee-\al_{+-abcd}^\vee-\al_{+-a'b'c'd'}^\vee\nn
 (a,b,c,d)&=&-(a',b',c',d')
\eea
Here there is an even number of minus signs amongst the signs $a,b,c,d$, and
therefore also amongst the signs $a',b',c',d'$. The result in Table 2 
is based on the choice $a=b=c=d=+$. More generally, the commutativity of
the translation operators immediately spoils an a priori possible uniqueness
of the procedure.
Furthermore, a very simple alternative to the general procedure for 
obtaining $t_{\al_i^\vee}$ exists when $\al_i\cdot \theta^\vee=1$, 
since in that
case the lemma applies directly on $\al_i$. In the following Table 3 we have 
listed those simple roots.\\[.2cm]
{\bf Table 3}
\begin{center}
\begin{tabular}{|c||c|c|c|c|c|c|c|c|c|}\hline
 &&&&&&&&&\\
 {\bf g}&\ $A_r$\ &\ $B_r$\ &\ $C_r$\ &\ $D_r$\ &\ $E_6$\ &\ $E_7$\ &\ 
  $E_8$\ &\ $F_4$\ &\ $G_2$\ \\
 &&&&&&&&&\\ \hline
 &&&&&&&&&\\
 $\ \al_i\ :\ \ \al_i\cdot\theta^\vee=1\ $  
 &
  $\al_1,\ \al_r$&$\al_2$&$\al_1$&$\al_2$&$\al_6$&$\al_6$&$\al_6$&$\al_2$
  &$\al_2$\\
 &&&&&&&&&\\ \hline
\end{tabular}\\[.5cm]
\end{center}
It will depend on the application which representation of the corresponding
translation operator $t_{\al_i^\vee}$ that is superior. Let us finally comment
on $G_2$ where the procedure does not apply directly due to
\ben
 t_{\theta^\vee}=t_{\al_1^\vee}\circ t_{\al_2^\vee}^2
\een
However, $\al_2\cdot\theta^\vee=1$ so we may simply use the lemma, the result
of which is in Table 2 of fundamental translations above. Similarly, 
the result for $t_{\al_1^\vee}$ in $C_r$ is obtained directly from the 
lemma, c.f. Table 3 above.

\subsection{Decompositions of Classical Weyl Reflections}

It remains to account for how to decompose the classical Weyl reflections.
For this purpose we use that if $\beta\cdot\al^\vee>0$ then 
$\beta\cdot\al^\vee=1,2,3$ and $\beta-\al\in\D$, 
and that from (\ref{easy}) it follows that
\ben
 \sigma_\beta=\sigma_\al\circ\sigma_{\beta-\beta\cdot\al^\vee\al}\circ
  \sigma_\al
\label{easy2}
\een
Let us illustrate our procedure for decomposing an arbitrary Weyl reflection
by considering the case of $B_r$ where the root type $\al_{ij}^-$ may be
expanded as
\ben
 \al_{ij}^-=\al_{j-1}+...+\al_i
\een
We see that
\ben
 (\al_{j-1}+...+\al_{i'})\cdot\al_{i'}^\vee=1\spa i'=i,...,j-2
\een
so according to (\ref{easy2}) the decomposition in Table 4 immediately 
follows. Similarly, the expansion 
\ben
 \al_{ii}=\al_r+...+\al_i
\een
leads to the decomposition of $\sigma_{\al_{ii}}$. Finally, we have
\ben
 \al_{ij}^+=\al_{ii}+\al_{jj}\spa\al_{ij}^+\cdot\al_{jj}^\vee=2
\een
so according to (\ref{easy2})
\ben
 \sigma_{\al_{ij}^+}=\sigma_{\al_{jj}}\circ\sigma_{\al_{ij}^+-2\al_{jj}}
  \circ\sigma_{\al_{jj}}=\sigma_{\al_{jj}}\circ\sigma_{\al_{ij}^-}\circ
  \sigma_{\al_{jj}}
\een
In Table 4, Table 5 and Table 6 we have summarized our findings for all
simple finite dimensional Lie algebras. Note that for notational reasons,
the decomposition of a Weyl reflection $\sigma_\al$ may depend on other
decompositions given (above it) in the tables.\\[.2cm]
\noindent {\bf Table 4}
\begin{center}
\begin{tabular}{|c||lr|}\hline
 &&\\
 {\bf g}&Decompositions of classical Weyl reflections  &\\ &&\\ \hline\hline
 &&\\
 $A_r$ &
  $\sigma_{\al_{ij}}=\sigma_i\circ...\circ\sigma_{j-2}\circ
  \sigma_{j-1}\circ\sigma_{j-2}\circ...\circ\sigma_i$ &\\
 &&\\ \hline
 &&\\
  $B_r$ &$\sigma_{\al_{ij}^-}=\sigma_i\circ...\circ\sigma_{j-2}\circ
  \sigma_{j-1}\circ\sigma_{j-2}\circ...\circ\sigma_i$  &\\ 
 & $\sigma_{\al_{ii}}=\sigma_i\circ...\circ\sigma_{r-1}\circ
  \sigma_{r}\circ\sigma_{r-1}\circ...\circ\sigma_i$  &$i\neq r$\\
 & $\sigma_{\al_{ij}^+}=\sigma_{\al_{jj}}\circ\sigma_{\al_{ij}^-}\circ
  \sigma_{\al_{jj}}$  &  \\
 &&\\ \hline &&\\
  $C_r$ &$\sigma_{\al_{ij}^-}=\sigma_i\circ...\circ\sigma_{j-2}\circ
  \sigma_{j-1}\circ\sigma_{j-2}\circ...\circ\sigma_i$ &\\
 &$\sigma_{\al_{ir}^+}=\sigma_{r}\circ...\circ\sigma_i\circ...\circ\sigma_r
  $&\\
 & $\sigma_{\al_{ij}^+}=\sigma_j\circ...\circ
  \sigma_{r-1}\circ\sigma_{\al_{ir}^+}\circ\sigma_{r-1}\circ
  ...\circ\sigma_j$&$j<r$\\
 & $\sigma_{\al_{ii}}=\sigma_{\al_{ir}^-}\circ\sigma_r\circ\sigma_{\al_{ir}^-}
  $&$i\neq r$\\
 &&\\ \hline
 &&\\
 $D_r$&$\sigma_{\al_{ij}^-}=\sigma_i\circ...\circ\sigma_{j-2}\circ
  \sigma_{j-1}\circ\sigma_{j-2}\circ...\circ\sigma_i$  &\\ 
 & $\sigma_{\al_{ir}^+}=\sigma_r\circ\sigma_{r-2}\circ...\circ\sigma_i
  \circ...\circ\sigma_{r-2}\circ\sigma_r$&$i<r-1$\\
 & $\sigma_{\al_{i,r-1}^+}=\sigma_r\circ\sigma_{r-1}\circ...\circ\sigma_i
  \circ...\circ\sigma_{r-1}\circ\sigma_r$&\\
 & $\sigma_{\al_{ij}^+}=\sigma_j\circ...\circ
  \sigma_{r-2}\circ\sigma_{\al_{i,r-1}^+}
  \circ\sigma_{r-2}\circ...\circ\sigma_j$ &$j<r-1$\\
  &&\\ \hline
\end{tabular}
\end{center}
\newpage
\noindent {\bf Table 5}
\begin{center}
\begin{tabular}{|c||ll|}\hline
 &&\\
 {\bf g}&Decompositions of classical Weyl reflections  &\\ &&\\ \hline\hline
 &&\\
  $E_6$ &$\sigma_{\al_{ij}^-}=\sigma_i\circ...\circ\sigma_{j-2}\circ
  \sigma_{j-1}\circ\sigma_{j-2}\circ...\circ\sigma_i   $  &\\ 
 & $\sigma_{\al_{1j}^+}=\sigma_6\circ\sigma_{\al_{2j}^-}\circ\sigma_6  $
  &$3\leq j\leq5$\\
 & $\sigma_{\al_{2j}^+}=\sigma_6\circ\sigma_{\al_{1j}^-}\circ\sigma_6  $ &\\
 & $\sigma_{\al_{3j}^+}=\sigma_2\circ\sigma_{\al_{2j}^+}\circ\sigma_2  $ &\\
 & $\sigma_{\al_{45}^+}=\sigma_3\circ\sigma_{\al_{35}^+}\circ\sigma_3  $ &\\
 & $\sigma_{\al_{\pm\pm\pm\pm\pm}}\mid_{(\mbox{{\footnotesize one}}\ +,\ 
  \mbox{{\footnotesize position}}\ i)}
  =\sigma_{\al_{1i}^-}\circ\sigma_5\circ\sigma_{\al_{1i}^-}$&$i\neq1$\\
 & $\sigma_{\al_{\pm\pm\pm\pm\pm}}\mid_{(\mbox{{\footnotesize three}}\ +,\ 
  \mbox{{\footnotesize positions}}
  \ i<j<k)}=\sigma_{\al_{jk}^+}\circ\sigma_5\circ\sigma_{\al_{jk}^+}$&$i=1$\\
 &$\sigma_{\al_{\pm\pm\pm\pm\pm}}\mid_{(\mbox{{\footnotesize three}}\ +,\ 
  \mbox{{\footnotesize positions}}
  \ i<j<k)}=\sigma_{\al_{jk}^+}\circ\sigma_{\al_{1i}^-}\circ\sigma_5\circ
  \sigma_{\al_{1i}^-}\circ\sigma_{\al_{jk}^+}$&$i\neq1$\\
 &$\sigma_{\al_{+++++}}=\sigma_{\al_{45}^+}\circ\sigma_{\al_{23}^+}
  \circ\sigma_5\circ\sigma_{\al_{23}^+}\circ\sigma_{\al_{45}^+}$&\\
 &&\\ \hline
 &&\\
 $E_7$ &
  $\sigma_{\al_{ij}^-}=\sigma_i\circ...\circ\sigma_{j-2}\circ
  \sigma_{j-1}\circ\sigma_{j-2}\circ...\circ\sigma_i  $ &\\
 & $\sigma_{\al_{1j}^+}=\sigma_7\circ\sigma_{\al_{2j}^-}\circ\sigma_7  $
  &$3\leq j\leq6$\\
 & $\sigma_{\al_{2j}^+}=\sigma_7\circ\sigma_{\al_{1j}^-}\circ\sigma_7  $ &\\
 & $\sigma_{\al_{3j}^+}=\sigma_2\circ\sigma_{\al_{2j}^+}\circ\sigma_2  $ &\\
 & $\sigma_{\al_{4j}^+}=\sigma_3\circ\sigma_{\al_{3j}^+}\circ\sigma_3  $ &\\
 & $\sigma_{\al_{56}^+}=\sigma_4\circ\sigma_{\al_{46}^+}\circ\sigma_4  $ &\\
 & $\sigma_{\al_{\pm\pm\pm\pm\pm\pm}}\mid_{(\mbox{{\footnotesize one}}\ +,\ 
  \mbox{{\footnotesize position}}\ i)}
  =\sigma_{\al_{1i}^-}\circ\sigma_6\circ\sigma_{\al_{1i}^-}$&$i\neq1$\\
 & $\sigma_{\al_{\pm\pm\pm\pm\pm\pm}}\mid_{(\mbox{{\footnotesize three}}\ +,\ 
  \mbox{{\footnotesize positions}}
  \ i<j<k)}=\sigma_{\al_{jk}^+}\circ\sigma_6\circ\sigma_{\al_{jk}^+}$&$i=1$\\
 &$\sigma_{\al_{\pm\pm\pm\pm\pm\pm}}\mid_{(\mbox{{\footnotesize three}}\ +,\ 
  \mbox{{\footnotesize positions}}
  \ i<j<k)}=\sigma_{\al_{jk}^+}\circ\sigma_{\al_{1i}^-}\circ\sigma_6\circ
  \sigma_{\al_{1i}^-}\circ\sigma_{\al_{jk}^+}$&$i\neq1$\\
 &$\sigma_{\al_{\pm\pm\pm\pm\pm\pm}}\mid_{(\mbox{{\footnotesize five}}\ +,\ 
  \mbox{{\footnotesize positions}}
  \ i<j<k<l<m)}=\sigma_{\al_{lm}^+}\circ
  \sigma_{\al_{jk}^+}\circ\sigma_6\circ
  \sigma_{\al_{jk}^+}\circ\sigma_{\al_{lm}^+}$ &$i=1$\\
 &$\sigma_{\al_{-+++++}}=\sigma_{\al_{56}^+}\circ
  \sigma_{\al_{34}^+}\circ\sigma_{\al_{12}^-}\circ\sigma_6\circ
  \sigma_{\al_{12}^-}\circ\sigma_{\al_{34}^+}\circ\sigma_{\al_{56}^+}$&\\
 &$\sigma_\theta=\sigma_6\circ\sigma_{\al_{-+++++}}\circ\sigma_6$&\\
 &&\\ \hline
\end{tabular}
\end{center}
\newpage 
\noindent {\bf Table 6}
\begin{center}
\begin{tabular}{|c||ll|}\hline
 &&\\
 {\bf g}&Decompositions of classical Weyl reflections  &\\ &&\\ \hline\hline
 &&\\
  $E_8$ &$ \sigma_{\al_{ij}^-}=\sigma_i\circ...\circ\sigma_{j-2}\circ
  \sigma_{j-1}\circ\sigma_{j-2}\circ...\circ\sigma_i $  &$j\neq8$\\ 
 & $\sigma_{\al_{1j}^+}=\sigma_8\circ\sigma_{\al_{2j}^-}\circ\sigma_8  $
  &$3\leq j\leq7$\\
 & $ \sigma_{\al_{2j}^+}=\sigma_8\circ\sigma_{\al_{1j}^-}\circ\sigma_8  $
  &$3\leq j\leq7$   \\
 &$\sigma_{\al_{ij}^+}=\sigma_{i-1}\circ\sigma_{\al_{i-1,j}^+}\circ
  \sigma_{i-1}$&$3\leq i\leq6$\\
 &&$j\neq8$\\
 & $\sigma_{\al_{\pm\pm\pm\pm\pm\pm\pm}}\mid_{(\mbox{{\footnotesize one}}\ +,\ 
  \mbox{{\footnotesize position}}\ i)}
  =\sigma_{\al_{1i}^-}\circ\sigma_7\circ\sigma_{\al_{1i}^-}$&$i\neq1$\\
 & $\sigma_{\al_{\pm\pm\pm\pm\pm\pm\pm}}
  \mid_{(\mbox{{\footnotesize three}}\ +,\ 
  \mbox{{\footnotesize positions}}
  \ i<j<k)}=\sigma_{\al_{jk}^+}\circ\sigma_7\circ\sigma_{\al_{jk}^+}$&$i=1$\\
 &$\sigma_{\al_{\pm\pm\pm\pm\pm\pm\pm}}
  \mid_{(\mbox{{\footnotesize three}}\ +,\ 
  \mbox{{\footnotesize positions}}
  \ i<j<k)}=\sigma_{\al_{jk}^+}\circ\sigma_{\al_{1i}^-}\circ\sigma_7\circ
  \sigma_{\al_{1i}^-}\circ\sigma_{\al_{jk}^+}$&$i\neq1$\\
 &$\sigma_{\al_{\pm\pm\pm\pm\pm\pm\pm}}\mid_{(\mbox{{\footnotesize five}}\ +,\ 
  \mbox{{\footnotesize positions}}
  \ i<j<k<l<m)}=\sigma_{\al_{lm}^+}\circ
  \sigma_{\al_{jk}^+}\circ\sigma_7\circ
  \sigma_{\al_{jk}^+}\circ\sigma_{\al_{lm}^+}$ &$i=1$\\
 &$\sigma_{\al_{\pm\pm\pm\pm\pm\pm\pm}}\mid_{(\mbox{{\footnotesize five}}\ +,\ 
  \mbox{{\footnotesize positions}}
  \ i<j<k<l<m)}$&\\
 &$\ \ \ =\sigma_{\al_{lm}^+}\circ
  \sigma_{\al_{jk}^+}\circ\sigma_{\al_{1i}^-}\circ\sigma_7\circ
  \sigma_{\al_{1i}^-}\circ
  \sigma_{\al_{jk}^+}\circ\sigma_{\al_{lm}^+}$ &$i\neq1$\\
 &$\sigma_{\al_{+++++++}}=\sigma_{\al_{67}^+}\circ
  \sigma_{\al_{45}^+}\circ\sigma_{\al_{23}^+}\circ\sigma_7\circ
  \sigma_{\al_{23}^+}\circ
  \sigma_{\al_{45}^+}\circ\sigma_{\al_{67}^+}$ &\\
 &$\sigma_{\al_{78}^-}=\sigma_7\circ\sigma_{\al_{-+++++-}}\circ\sigma_7$&\\
 &$\sigma_{\al_{i8}^-}=\sigma_i\circ\sigma_{\al_{i+1,8}^-}\circ\sigma_i$&
  $i<7$\\
 &$\sigma_{\al_{18}^+}=\sigma_7\circ\sigma_{\al_{+++++++}}\circ\sigma_7$&\\
 &$\sigma_{\al_{i8}^+}=\sigma_{i-1}\circ\sigma_{\al_{i-1,8}^+}\circ
  \sigma_{i-1}$&$2\leq i\leq7$\\
 &&\\ \hline 
 &&\\
 $F_4$ &$\sigma_{\al_{--+}}=\sigma_4\circ\sigma_1\circ\sigma_4$ 
  \hspace{2.3cm}$
  \sigma_{\al_{-+-}}=\sigma_3\circ\sigma_{\al_{--+}}\circ\sigma_3 $ &\\
 &$ \sigma_{\al_{-++}}=\sigma_4\circ\sigma_{\al_{-+-}}\circ\sigma_4$ 
  \hspace{1.7cm}$
  \sigma_{\al_{+-+}}=\sigma_2\circ\sigma_{\al_{-++}}\circ\sigma_2 $&\\
 &$ \sigma_{\al_{++-}}=\sigma_3\circ\sigma_{\al_{+-+}}\circ\sigma_3$
  \hspace{1.7cm}$
  \sigma_{\al_{+++}}=\sigma_4\circ\sigma_{\al_{++-}}\circ\sigma_4   $&\\
 & $\sigma_{\al_{+--}}=\sigma_2\circ\sigma_{\al_{-+-}}\circ\sigma_2   $&\\
 & $\sigma_{\al_{11}}=\sigma_1\circ\sigma_{\al_{+++}}\circ\sigma_1$
  \hspace{2.0cm}$
 \sigma_{\al_{22}}=\sigma_{\al_{-++}}\circ\sigma_{\al_{+++}}\circ
  \sigma_{\al_{-++}} $&\\ 
 &$ \sigma_{\al_{33}}=\sigma_{\al_{+-+}}\circ\sigma_{\al_{+++}}\circ
  \sigma_{\al_{+-+}} $&\\
 &$ \sigma_{\al_{24}^-}=\sigma_3\circ\sigma_2\circ\sigma_3 $
  \hspace{2.6cm}$
  \sigma_{\al_{12}^-}=\sigma_{\al_{--+}}\circ\sigma_3\circ
  \sigma_{\al_{--+}} $&\\
 &$ \sigma_{\al_{13}^-}=\sigma_{\al_{--+}}\circ\sigma_{\al_{24}^-}\circ
  \sigma_{\al_{--+}}$
  \hspace{1cm}$
 \sigma_{\al_{14}^-}=\sigma_{\al_{-+-}}\circ\sigma_2\circ
  \sigma_{\al_{-+-}}$ &\\
 & $\sigma_{\al_{ij}^+}=\sigma_{\al_{jj}}\circ\sigma_{\al_{ij}^-}\circ
  \sigma_{\al_{jj}}$&\\
 &&\\ \hline
 &&\\
  $G_2$ &$\sigma_{\al_{13}}=\sigma_2\circ\sigma_1\circ\sigma_2$
  \hspace{2.7cm}$
  \sigma_{\al_{23}}=\sigma_1\circ\sigma_2\circ\sigma_1\circ\sigma_2
  \circ\sigma_1  $&  \\ 
 & $\sigma_{\al_{+-+}}=\sigma_1\circ\sigma_2\circ\sigma_1$
  \hspace{2.4cm}$
  \sigma_\theta=\sigma_2\circ\sigma_1\circ\sigma_2\circ\sigma_1\circ
  \sigma_2  $ &\\
 &&\\ \hline
\end{tabular}
\end{center}
\newpage
\noindent Obviously, these decompositions are not unique. 
In Table 4, Table 5 and Table 6 we have only listed one solution for each 
type of roots and not included all
the fundamental Weyl reflections themselves.

\section{Conclusions}
In this paper we have provided explicit
decompositions of affine Weyl reflections
in terms of fundamental reflections. Together with the differential operator
realizations of the Lie algebras \cite{Ras1,PRY}, this result allows 
generalizing the work \cite{AY} on fusion rules in affine $SL(2)$ current
algebra to more general affine current algebras.
We intend to come back elsewhere with a discussion along these lines
\cite{Ras3}. 

Generalizations to superalgebras are also interesting. Recently, we
have presented differential operator realizations of basic Lie superalgebras
\cite{Ras2}. However, in order to determine the fusion rules in affine
current superalgebras along the lines of Awata and Yamada \cite{AY} we
need to work out
not only decompositions of affine Weyl (super-)reflections but also
Malikov-Feigin-Fuks type \cite{MFF} expressions for the singular vectors. 
Studies of these issues would be interesting to undertake in future works.
\\[.5 cm] 
{\bf Acknowledgement}\\[.2cm]
The author thanks V.K. Dobrev for a copy of the unpublished preprint
\cite{Dob} and for pointing out reference \cite{DS}.
\newpage
\appendix
\section{Embedding of Root Systems}

In this appendix we review the well known embeddings (see e.g. \cite{Fuc})
of the root systems in $\R^n$ by means of an orthonormal basis
$\{e_i\}_{i=1,...,n}$. In the case of ${\bf \mbox{g}}=B_r,C_r,D_r,E_8,F_4$
we have $n=r$, in the case of ${\bf \mbox{g}}=A_r,E_7,G_2$ we have $n=r+1$
while for ${\bf \mbox{g}}=E_6$ $n=r+2$. The normalizations of the root
systems are such that the short roots have length squared equal to 2,
except for $B_r$ and $F_4$ where the long roots have length squared
equal to 2.
However, a fixing of the root lengths here is of no importance
for the general arguments in the main body of this paper where we keep
the normalization of the root system ($\theta^2$) a free parameter.
The explicit embeddings are given in the following tables, Table 7,
Table 8 and Table 9.\\[.2cm]
\noindent {\bf Table 7}
\begin{center}
\begin{tabular}{|c||lr|l|}  \hline
 &&&\\
 {\bf g} &  Simple roots & & $\D_+$ \\
 &&&\\ \hline\hline
 &&&\\
 $A_r$   &  $\al_i=e_i-e_{i+1}$  &  $1\leq i\leq r$  
         &  $\al_{ij}=e_i-e_j$\hspace{.9cm}  $1\leq i<j\leq r+1$\\
 &&&\\ \cline{4-4}
 &&&\\
         &                       &
         &  $\theta=\sum_{i=1}^r\al_i=e_1-e_{r+1}$    \\
 &&&\\ \hline
 &&&\\
 $B_r$   &  $\al_i=e_i-e_{i+1}$           &  $1\leq i< r$
         &  $\al_{ij}^{\pm}=e_i\pm e_j$   \hspace{1cm} $1\leq i<j\leq r$ \\
         &  $\al_r=e_r$          &
         &  $\al_{ii}=e_i$   \hspace{2.3cm} $1\leq i\leq r$    \\
 &&&\\ \cline{4-4}
 &&&\\
         &                       &
         &  $\theta=\al_1+2\sum_{i=2}^r\al_i=e_1+e_2$    \\
 &&&\\ \hline
 &&&\\
 $C_r$   &  $\al_i=e_i-e_{i+1}$  &  $1\leq i<r$
         &  $\al_{ij}^{\pm}=e_i\pm e_j$ \hspace{1cm}  $1\leq i<j\leq r$ \\
         &  $\al_r=2e_r$          &
         &  $\al_{ii}=2e_i$   \hspace{2.1cm} $1\leq i\leq r$    \\
 &&&\\ \cline{4-4}
 &&&\\
         &                       &
         &  $\theta=2\sum_{i=1}^{r-1}\al_i+\al_r=2e_1$       \\
 &&&\\ \hline
 &&&\\ 
 $D_r$   &  $\al_i=e_i-e_{i+1}$  &  $1\leq i< r$
         &  $\al_{ij}^{\pm}=e_i\pm e_j$  \hspace{1cm}  $1\leq i<j\leq r$ \\
         &  $\al_r=e_{r-1}+e_r$           &
         &                                            \\
 &&&\\ \cline{4-4}
 &&&\\
         &                       &
         &  $\theta=\al_1+2\sum_{i=2}^{r-2}\al_i+\al_{r-1}+\al_r$\\
 &&&\hspace{.2cm}                   $=e_1+e_2$                    \\
 &&&\\ \hline
\end{tabular}
\end{center}
\newpage
\noindent {\bf Table 8}
\begin{center}
\begin{tabular}{|c||lr|l|}  \hline
 &&&\\
 {\bf g} &  Simple roots & & $\D_+$ \\
 &&&\\ \hline\hline
 &&&\\
 $E_6$   &  $\al_i=-e_i+e_{i+1}$   &  $1\leq i\leq 4$
         &  $\al_{ij}^{\pm}=\pm e_i+e_j$ \hspace{1cm}  $1\leq i<j\leq 5$ \\
         &  $\al_5=\hf\left(e_1-\sum_{j=2}^7e_j+e_8\right)$          &
         &  $\al_{\pm\pm\pm\pm\pm}=
            \hf\left(\sum_{j=1}^5(\pm)e_j-e_6-e_7+e_8\right)$  \\
            &&&   even number of minus signs   \\
         &  $\al_6=e_1+e_2$          &
         &                   \\
 &&&\\ \cline{4-4}
 &&&\\
         &                       &
         &  $\theta=2\al_1+3\al_2+2\al_3$\\
 &&&\hspace{.3cm}$+\ \al_4+\al_5+2\al_6$       \\
 &&&       \hspace{.3cm}$=\hf\left(\sum_{j=1}^5e_j-e_6-e_7+e_8\right)$\\
 &&&\\ \hline
 &&&\\
 $E_7$   &  $\al_i=-e_i+e_{i+1}$   &  $1\leq i\leq 5$
         &  $\al_{ij}^{\pm}=\pm e_i+e_j$ \hspace{1cm}  $1\leq i<j\leq 6$ \\
         &  $\al_6=\hf\left(e_1-\sum_{j=2}^7e_j+e_8\right)$          &
         &  $\al_{78}=-e_7+e_8$                 \\
         &  $\al_7=e_1+e_2$          &
         &  $\al_{\pm\pm\pm\pm\pm\pm}=
                  \hf\left(\sum_{j=1}^6(\pm)e_j-e_7+e_8\right)$  \\
            &&&   odd number of minus signs   \\
  &&&\\ \cline{4-4}
 &&&\\
         &                       &
         &  $\theta=3\al_1+4\al_2+3\al_3$\\
 &&&\hspace{.3cm}$+\ 2\al_4+\al_5+2\al_6+2\al_7$       \\
 &&&       \hspace{.3cm}$=-e_7+e_8$\\
 &&&\\ \hline
 &&&\\
 $E_8$   &  $\al_i=-e_i+e_{i+1}$   &  $1\leq i\leq 6$
         &  $\al_{ij}^{\pm}=\pm e_i+e_j$ \hspace{1cm}  $1\leq i<j\leq 8$ \\
         &  $\al_7=\hf\left(e_1-\sum_{j=2}^7e_j+e_8\right)$          &
         &  $\al_{\pm\pm\pm\pm\pm\pm\pm}=
             \hf\left(\sum_{j=1}^7(\pm)e_j+e_8\right)$  \\
            &&&   even number of minus signs   \\
         &  $\al_8=e_1+e_2$          &
         &                   \\
 &&&\\ \cline{4-4}
 &&&\\
         &                       &
         &  $\theta=4\al_1+6\al_2+5\al_3+4\al_4$\\
 &&&\hspace{.3cm}$+\ 3\al_5+2\al_6+2\al_7+3\al_8$       \\
 &&&       \hspace{.3cm}$=e_7+e_8$\\
 &&&\\ \hline
\end{tabular}
\end{center}
\newpage
\noindent {\bf Table 9}
\begin{center}
\begin{tabular}{|c||lr|l|}  \hline
 &&&\\
 {\bf g} &  Simple roots & & $\D_+$ \\
 &&&\\ \hline\hline
 &&&\\
 $F_4$   &  $\al_1=\hf\left(e_1-e_2-e_3-e_4\right)$   & 
         &  $\al_{\pm\pm\pm}=\hf\left(e_1\pm e_2\pm e_3\pm e_4\right)$ \\
         &  $\al_2=e_2-e_3$          &
         & $\al_{ij}^\pm=e_i\pm e_j$\hspace{1.3cm}$1\leq i<j\leq4$  \\
         &  $\al_3=e_3-e_4$          &
         & $\al_{ii}=e_i$\hspace{2.6cm}$1\leq i\leq4$                 \\
         &  $\al_4=e_4$  &&\\
 &&&\\ \cline{4-4}
 &&&\\
         &                       &
         &  $\theta=2\al_1+2\al_2+3\al_3+4\al_4$\\
 &&&       \hspace{.3cm}$=e_1+e_2$\\
 &&&\\ \hline
 &&&\\
 $G_2$   &  $\al_1=e_1-e_2$   & 
         &  $\al_1,\ \al_2$ \\
         &  $\al_2=-2e_1+e_2+e_3$          &
         & $\al_{13}=-e_1+e_3,\ \al_{23}=-e_2+e_3$  \\
         &            &
         & $\al_{+-+}=e_1-2e_2+e_3$                 \\
         & && $\al_{--+}=-e_1-e_2+2e_3$\\
 &&&\\ \cline{4-4}
 &&&\\
         &                       &
         &  $\theta=3\al_1+2\al_2=-e_1-e_2+2e_3$\\
 &&&\\ \hline
\end{tabular}
\end{center}

\end{document}